\documentclass[12pt, letterpaper]{article}

\usepackage{amssymb}
\usepackage{amsfonts, amsmath}

\begin{document}
\begin{center}
\textbf{The Sums of the $k$-powers of the Euler set and their connection with Artin's conjecture for primitive roots}\\\vspace{0.5cm}
Constantin M. Petridi
\\ cpetridi@math.uoa.gr
\\ cpetridi@hotmail.com
\end{center}\vspace{0.5cm}

\begin{flushleft}
$``$Deepest interrelationships in analysis are of an arithmetical nature".\\
Hermann Minkowski, preface to $``$Diophantiche Approximationen".
\end{flushleft}\vspace{1cm}
\textbf{Abstract}. We examine the sums $\mathcal{S}(k,\,n)$ of the $k$-th powers of the $\phi(n)$ integers $\alpha_{1}<\alpha_{2}<\cdots<\alpha_{\phi(n)}$ less than and prime to $n$ (Euler set [7]) and prove a formula (new) for $\mathcal{S}(3,\,n)$. If $n$ equals a prime $p$, we prove a theorem showing a connection of $\mathcal{S}(3,\,p)$ with Artin's conjectural constant for primitive roots and with other functions involving primes.\vspace{1cm}\\

\hspace{-0.7cm}\textbf{1 Introduction}\vspace{0.5cm}\\
Let $\phi(n)$ be the Euler totient function, $\alpha_{1}(=1)<\alpha_{2}<\cdots<\alpha_{\phi(n)}(=\alpha-1)$ the $\phi(n)$ integers less than and prime to $n=p_{1}^{e_{1}}p_{2}^{e_{2}}\cdots p_{w}^{e_{w}}$, $R(n)=p_{1}p_{2}\cdots p_{w}$ the radikal of $n$. We examine the sums
\[ \mathcal{S}(k,\,n)=\sum_{i=1}^{\phi(n)}\alpha_{i}^{k},\,\,k=0,\,1,\,2,\cdots.\]
Known are following formulas\vspace{0.5cm}
\begin{eqnarray*}
\mathcal{S}(0,\,n)&=&\phi(n)\\
\mathcal{S}(1,\,n)&=&\frac{1}{2}\phi(n)n\\
\mathcal{S}(2,\,n)&=&\frac{1}{3}\phi(n)\big(n^{2}+{(-1)}^{w}\frac{1}{2}R(n)\big).
\end{eqnarray*}\\
\vspace{0.5cm}The first two are trivial. The third was proved by G. {P$\acute{o}$lya}\footnote{Ronald C. Read in his preface to $[6]$, writes
$``$ In 1937, there appeared a paper that was to have a profound influence on the progress of combinatorial enumeration. Entitled Kombinatorische Anzahlbestimmungen f$\ddot{u}$r Gruppen, Graphen und Chemische Verbindungen. It was published in Acta Mathematica, Vol. 62, pp. 145 to 254. Its author George P$\acute{o}$lya ...$"$. P$\acute{o}$lya's paper gave rise to the innovative researches of Rota, Sch$\ddot{u}$tzenberger, Wilf, Petrov$\check{s}$k, Zeilberger, Stanley, Bj$\ddot{o}$rner et al.}
and G. Szeg$\ddot{o}$ in their book $[5]$.\\
We here prove that\\
\[ \mathcal{S}(3,\,n)\,=\,\frac{1}{4}\phi(n)\big(n^{3}+(-1)^{w}R(n)n\big), \]\\
and show that there is a connection between $\mathcal{S}(3,\,p)$, $p$ prime, and Artin's conjecture for primitive roots. In conclusion, we mention, without discussing, that analogous connections exist for other conjectures resp. functions involving primes, f.ex. $\zeta(s)$.\vspace{1cm}\\

\hspace{-0.7cm}\textbf{2 Formula for  $\mathcal{S}(3,\,n)$ and conjecture for $\mathcal{S}(k,\,n)$}\vspace{0.5cm}\\
\textbf{Theorem 1}\vspace{0.5cm}\\
Let $\{\alpha_{1},\,\alpha_{2},\,\cdots,\,\alpha_{\phi(n)}\}$ be the Euler set of $n$. Then\\
\[ \mathcal{S}(3,\,n)\,=\,\frac{1}{4}\phi(n)\big(n^{3}+(-1)^{w}R(n)n\big). \]\vspace{0.1cm}\\
\hspace{-0.7cm}\textbf{Proof}. Literally copying from $[5]$, (Number Theory, CH. VIII, exercise 27) we define
\[
\begin{array}{ccc}
\text{Objects\,\,are}\,\,\,\,\,\,: & & 1,\,2,\,\cdot\cdot\cdot,m,\,\,m\in\mathbb{Z}\\
\\
\text{Properties\,\,are}: & & \text{divisibility\,\,by}\,\,p_{i},\,i=1,\,2,\,\cdots,\,w\\
\\
\hspace{-1cm}\text{Value}\,\,W(m)\,\,\text{of}\,\,W\,\text{is}: & & W(m)=1^3+2^3+\cdot\cdot\cdot+m^3=\dfrac{(m(m+1))^2}{4}\\ \\
                                                       & & \hspace{1.2cm}=Am^4+B m^3+C m^{2}+D m+E\\ \\
                                                       & & \hspace{0.8cm}A=B=C=\frac{1}{4},\,\,\,\,\,\,\,D=E=0.
\end{array}
\]
Inserting for $m$ the primes $p_1,\,p_2,\,\cdots,p_w$ in the values $W(m)$ we have\vspace{0.5cm}
\begin{eqnarray*}
W(p_1)&=&Ap_{1}^{4}+Bp_{1}^{3}+C p_{1}^{2}=\frac{1}{4}(p_{1}^{4}+p_{1}^{3}+p_{1}^{2})\\
W(p_2)&=&Ap_{2}^{4}+Bp_{2}^{3}+C p_{2}^{2}=\frac{1}{4}(p_{2}^{4}+p_{2}^{3}+p_{2}^{2})\\
\cdots& &\cdots\cdots\cdots\cdots\cdots\cdots\cdots\cdots\cdots\cdots\cdots\\
W(p_w)&=&Ap_{w}^{4}+Bp_{w}^{3}+C p_{w}^{2}=\frac{1}{4}(p_{w}^{4}+p_{w}^{3}+p_{w}^{2}).
\end{eqnarray*}\\
Application of the inclusion-exclusion, as done in $[5]$ by P$\acute{o}$lya and Szeg$\ddot{o}$ for the case $k=2$ gives
\[ \mathcal{S}(3,\,n)\,=\,\frac{1}{4}\phi(n)\big(n^{3}+(-1)^{w}R(n)n\big), \]
\\
using the facts that\vspace{0.5cm}
\begin{eqnarray*}
1-\sum_{i=1}^{\omega}\frac{1}{p_i}+\sum_{i,\,j=1}^{\omega}\frac{1}{p_{i}p_{j}}-\cdots&=&\frac{\phi(n)}{n}\\
1-\sum_{i=1}^{\omega}{p_i}+\sum_{i,\,j=1}^{\omega}{p_{i}p_{j}}-\cdots&=&\frac{\phi(n)}{n}(-1)^{\omega}R(n).
\end{eqnarray*}
Q.E.D.\\ \\ \\
For exponents $k$ higher than $3$, above method fails because in the inclusion-exclusion summation occur sums over the primes $p_1,\,p_2,\,\cdots,\,p_{w}$ which are not expressible in terms of $n,\,\phi(n),\,(-1)^w$ and $R(n)$, as is the case for $k\leq{3}$.\\ \\
Adding above formula $\mathcal{S}(3,\,n)$ to the known formulas we get following list\vspace{0.3cm}\\
\begin{eqnarray*}
\mathcal{S}(0,\,n)&=&\phi(n)\\ \\
\mathcal{S}(1,\,n)&=&\frac{1}{2}\phi(n)n\\ \\
\mathcal{S}(2,\,n)&=&\frac{1}{3}\phi(n)\big(n^2+(-1)^w\frac{1}{2}R(n)\big)\\ \\
\mathcal{S}(3,\,n)&=&\frac{1}{4}\phi(n)\big(n^3+(-1)^{w}R(n)n\big)
\end{eqnarray*}\\
Numerical computations results for $k\geq4$ lead us to believe that in general we have\\
\[
\mathcal{S}(k,\,n)\,=\,\frac{1}{k+1}\phi(n)\big(n^{k}+c_{1}n^{k-1}+\cdots+c_{k}\big),\vspace{0.5cm}
\]
where $c_1,\,c_2,\,\cdots,\,c_k$ are constants depending linearly on $k,\,R(n),\,(-1)^w$. If true this would give\vspace{0.5cm}\\
\[ \lim_{n\to\infty}\frac{S(k,\,n)}{\phi(n)n^{k}}\,=\,\frac{1}{k+1}.\]\vspace{1cm}\\

\hspace{-0.7cm}\textbf{3 Connection of $\mathcal{S}(3,\,p)$ with Artin's conjecture for primitive
roots [1], [2], [3], [4]}\vspace{0.5cm}\\
\textbf{Theorem 2}\vspace{0.5cm}\\
Denoting by $C$ the convergent (as seen below) product
\[C=\prod_{p\,\text{prime}}\big(1-\frac{1}{p(p-1)}\big)\] we have
\[
\prod_{p\,\text{prime}}\big(1-\frac{1}{2\sqrt{\mathcal{S}(3,\,p)}}\big)=C.
\]\\
\hspace{-0.5cm}\textbf{Proof}. In section $2$ we proved that
\[ \mathcal{S}(3,\,n)=\frac{1}{4}\phi(n)(n^3+(-1)^wR(n)n). \]
For $n$ equal to a prime $p$ this becomes
\[ \mathcal{S}(3,\,p)=\frac{1}{4}(p-1)(p^{3}-p^{2})=\frac{1}{4}(p(p-1))^2. \]
Or
\[ 2\sqrt{\mathcal{S}(3,\,p)}=p(p-1) \]
and hence
\[1-\frac{1}{2\sqrt{\mathcal{S}(3,\,p)}}=1-\frac{1}{p(p-1)}. \]
Taking the products over all primes of both sides gives\\ \\
\[ \prod_{p\,\text{prime}}\big(1-\frac{1}{2\sqrt{\mathcal{S}(3,\,p)}}\big)=\prod_{p\,\text{prime}}\big(1-\frac{1}{p(p-1)}\big)=C.\]\\ \\
But according to [4] $C$ is Artin's constant $C_{\text{Artin}}=0.3739558136\cdots$ (sequence $A005596$ in On-line Encyclopedia of Integer Sequences).\\ \\ \\
Q.E.D.\vspace{1cm}\\

\hspace{-0.7cm}\textbf{4 Connections of $\mathcal{S}(3,\,p)$ with functions involving primes}\vspace{0.5cm}\\
Pursuing the same line of thought as before but solving $\mathcal{S}(3,\,p)=\frac{1}{4}(p(p-1))^2$
for $p$ instead of for $p(p-1)$, we have\\ \\
\[p=\frac{1}{2}+\frac{1}{2}\sqrt{1+4\sqrt{\mathcal{S}(3,\,p)}}.\]\\ \\
This value of a prime in terms of $S(3,\,p)$ can now be substituted in functions involving the set of primes or its subsets. An immediate example at hand is Riemann's Zeta function\\
\[ \zeta(s)=\prod_{p\,\text{prime}}\frac{1}{1-p^{-s}}. \]\\ \\
Inserting for $p$ its expression in terms of $\mathcal{S}(3,\,p)$, we get\\ \\
\[  \zeta(s)=\prod_{p\,\text{prime}}\frac{1}{1-\big\{\frac{1}{2}+\frac{1}{2}\sqrt{1+4\sqrt{\mathcal{S}(3,\,p)}}\big\}^{-s}}. \]\\ \\
Clearly, above procedure can be applied to other functions resp. conjectures whose analytical expressions involve primes.\\ \\ \\
\textbf{Note}. What we said in the last two sections are in no way to be interpreted as meaning that our statements are of value for theoretical or computational purposes in the theory of prime numbers. This is an open question pending decision, if at all. In the meantime we would say that these statements just exemplify Henri Poincar$\acute{e}$s dictum $``$Mathematics is the art of giving the same name to different things$"$ (La Science et L'hypoth$\grave{e}$se, 1902).\vspace{0.5cm}\\

\hspace{-0.7cm}\textbf{References}\vspace{0.5cm}\\
$[1]$ {E. Artin, Collected papers, Addison Wesley, 1965, pp. VIII-X.}\\ \\
$[2]$ {M. Ram Murty, On Artin's Conjecture, Journal of Number Theory 16 (1983),\\ pp. 147-168.}\\ \\
$[3]$ {M. Ram Murty, Artin's Conjecture for Primitive Roots, Math. Intellingencer, 10 (1988)\\ pp. 59-67.}\\ \\
$[4]$ {\text{https:}\, $//$ \,\text{en.wikipedia.org/Artin,}\,$\%$\,\,{\text{27S\_conjecture}\_on\_primitive\_roots}.\\ \text{This page was last modified on 13 October 2016, at 04:07}.\\ \\
$[5]$ {G. P$\acute{o}$lya and G. Szeg$\ddot{o}$, Aufgaben und Lehrs$\ddot{a}$tze aus der Analysis, Berlin, Verlag von Julius Springer, Berlin, 1925.}\\ \\
$[6]$ {Ronald C. Read, Preface to G. P$\acute{o}$lya, \hspace{0.2cm}R. C. Read, Combinatorial Enumeration of Groups, Graphs and Chemical Compounds, 1987, Springer Verlag, New-York, Inc.}\\ \\
$[7]$ {Constantin M. Petridi, The integer recurrence $P(n)=a+P\big(n-\phi(a)\big)$ I, arXiv: 1208.5348v3 [math.NT] 4 Feb 2014.}

\end{document}